\documentclass[11pt]{article}
\usepackage[utf8]{inputenc}
\usepackage[
papersize={5.5in, 8.5in},
left=0.35in,
right=0.35in,
top=0.35in,
bottom=0.5in]{geometry}
\usepackage{graphicx}
\usepackage{amsthm}
\usepackage{rotating}
\usepackage{amsfonts,amssymb,amsmath}
\usepackage[usenames]{color}

\usepackage{array}
\usepackage{float}
\usepackage{booktabs}
\usepackage[dvipsnames,usenames,table,xcdraw]{xcolor}

\graphicspath{ {images/} }

\title{\bf 
What percent of the plane can be properly 5- and 6-colored?
}
\author{\bf 
\textcolor[rgb]{0,0.75,0}{Jaan Parts} \\
} 
\date{\normalsize \textcolor[rgb]{0,0.75,0}{Kazan, Russia, jaan\_parts@.mail.ru}}

\begin{document}

\maketitle

\pagestyle{empty}
\thispagestyle{empty}

\begin{abstract}
We present a tiling of more than 99.985698\% of the Euclidean plane with six colors, reducing the previous record for uncovered fraction of the plane by about 12.8\%. We also present a tiling of more than 95.99\% of the plane with five colors. It is thus shown that any unit-distance graph of order at most 6992 and 24 in the plane can be properly 6-colored and 5-colored, respectively.
\end{abstract}

\section{Preliminaries}
In this paper we will consider the following main problem: what maximum fraction of the Euclidean plane can be properly covered by tiles of a given number of colors?

It is known that 7 colors are enough to cover the entire plane (below we give some non-trivial examples). In the case of 6 or fewer colors, one may ask how much of the plane remains uncolored? One of the first such estimates was given by Edward Pegg, Jr. \cite[p.~18]{soi}, proposed a tiling based on heptagons. Dan Pritikin who gave a stronger estimate obtained on its basis a lower bound of the order of an arbitrary 6-colorable unit-distance graph \cite{pri}. Here we take a small step forward.

\paragraph{Definitions.}
A \textit{tile-based coloring}, or \textit{tiling}, of the plane is its partition by polygons or other connected geometric shapes, called \textit{tiles}, with no overlaps and no gaps, together with the assignment of a color to each tile.
A \textit{$k$-tiling} is one that uses exactly $k$ colors.

A \textit{proper tiling} of the plane is a tiling in which no two points of the plane at unit-distance from each other have the same color.

A \textit{partial tiling} is a proper tiling of (not necessarily connected) subset of the plane (rather than the whole plane). 
A \textit{void} is a connected region of the non-tiled part of the plane.

A pair of tiles or a pair of voids all points of which are located at most at a unit-distance is termed \textit{Siamese}.
A tiling is termed \textit{Siamese-free} if it features no Siamese pair of tiles.

A tile that lies on one side of any line containing it is termed \textit{convex}.

A tiling of the plane can be considered as an infinite planar graph with all vertices having degree at least 3. Here tile \textit{perimeter} corresponds to graph cycle that do not contain any smaller cycles. To avoid misinterpretation, we will explicitly distinguish the terminology associated with the descriptions of a tiling as a planar graph (with vertices and edges) and as a system of separate polygons (with corners and sides).
A \textit{vertex} is the connection point of three or more \textit{edges} of a graph.
A \textit{corner} is the junction point of two \textit{sides} of a polygon.
An edge may thus contain multiple corners. And a side of a tile may contain multiple vertices.

\paragraph{Designations.}
To compare the results, it is more convenient to use as a value being maximized the ratio $\rho=S_{\Sigma}/S_{\Delta}$, where $S_{\Delta}$ is the area of all voids of some repeating part of the tiling with total area $S_{\Sigma}$. Accordingly, the fraction of the plane that is occupied by voids is $\delta=1/\rho=S_{\Delta}/S_{\Sigma}$.

The number of available colors is denoted by $k$. The colors are numbered from 1 to $k$. 
And $\rho_k$ means $\rho$ for $k$ colors.

For the notation of base points (including corners, vertices and some others) it is more convenient to snap not to tiles, but to voids. Voids and individual points are denoted by capital letters in alphabetical order, starting with $A$. The first few letters correspond to different types of voids that are not equivalent up to translation and reflection. The points on the perimeter of a void are denoted by the void designation with an additional index, as a letter-number pair, and are numbered clockwise, starting from a specially highlighted point with number 1. If it is necessary to list several points of the same void sequentially, then the void is indicated once, and the necessary indices are separated by commas, for example, $A_1A_2A_3A_4$ is written $A_{1,2,3,4}$.

The repeating unit of a tiling (a \textit{cell}) is highlighted by dimming together with a dashed line. Here we do not pay attention to the tile colors and use tiling symmetries to reduce the size of the cell.

\paragraph{Constraints.}

As a rule, each tile tends to increase its area and come closer to another tile of the same color in order to maximize the fraction of the covered plane.
But tile growth is limited by tiles of other colors and \textit{distance constraints} arising from the unit-distance requirement. As result tile boundaries take on a complex shape.
We distinguish two main types of distance constraints: \textit{inner} (tile width less than 1) and \textit{outer} (distance between tiles greater than 1).
If two colors alternate along a straight line, then a special \textit{outer-pair} constraint appears: the distance between the edges with the same color transition is at least 2.
In the case of Siamese tiles, a special combination of inner and outer constraints arises: distance between any two tiles either less than 1 or greater than 1.

It is often convenient to consider the constraints in groups of 4 points presented as quadrilateral $ABCD$\footnote{The order of points is important; $AD$ and $BC$ intersect for inner constraints and do not intersect for outer ones.} that determine the distance between line segments or arcs $AB$ and $CD$. For inner constraints the distances $|AD|$ and $|BC|$ are at most 1; for outer ones $|AD|$ and $|BC|$ are at least 1 (or 2 for outer-pair constraints). 
Usually it is enough to consider the special case where $ABCD$ is an isosceles trapezoid, and use equalities rather than inequalities, assuming the corresponding curving and coloring of sides $AB$ and $CD$.
If only convex tiles are acceptable, then extended \textit{outer-ext} constraints are needed on the distance between non-parallel line segments. 
We also consider the cases $ABC = ABCC$ and $BC = BBCC$, where the trapezoid degenerates into a triangle or a line segment.
In the figures, spacing quadrilaterals are shown by colored line segments: inner ones are black, outer are red, outer-pair are blue. 

\paragraph{Curving.} \label{curving}
In the simplest case, tiles are a polygons with straight sides. More generally, the sides of a tile are \textit{arcs} of a certain radius. We will need to distinguish between two types of segment rounding: curved and wavy. 
A \textit{curved} side always separates a tile from a void, and has one optimal shape (is \textit{rigid} in this sense).
A \textit{wavy} edge always separates two tiles, and typically splits into the rigid segment and the remaining non-rigid parts, which have a range of shape options that do not affect the fraction of the covered plane (and are highlighted by dashed lines in the figures).

Wavy edges are formed under the influence of outer-pair constraints and typically lead to non-convex tiles. In the optimization procedure, the wavy edge is involved as a whole, without dividing into parts.

Curved sides also matter in the pursuit of a record, so we briefly consider the options.
Let a tile have two parallel sides $AB$ and $CD$, $AD$ and $BC$ are intersected, $|AD|=|BC|=1$, $|AC|=|BD|$, $|AB|=a$, $|CD|=b$. We compare two curving methods such that the distance between the curved sides does not exceed 1. The first way is to draw two arcs of radius 1 centered at points $A$ and $B$ from points $D$ and $C$ until they intersect at point $E$ from which an arc of radius 1 is drawn between the points $A$ and $B$. The second way is to find the intersection $O$ of the diagonals $AD$ and $BC$ and draw arcs with radii $r=|OA|$ and $(1-r)=|OD|$. (In passing, we note that we cannot complicate the construction and, for example, move $O$ when drawing arcs, since this will lead to an increase in the maximum distance between points on the arcs.) 

The increment of the tile area in the first ($S_1$) and second ($S_2$) method is described by the formulae
$$ S_1=\frac{b h}{2}+\beta–\sin\beta+\frac{\gamma–\sin\gamma}{2},\quad S_2=\left(\frac{1}{2}-\frac{a b}{(a+b)^2}\right)(\alpha-\sin\alpha),$$ 
where 
$h\!=\!\sqrt{1-\left(\frac{a}{2}\right)^2}-\sqrt{1-\left(\frac{a+b}{2}\right)^2}$,
$\alpha\!=\!\angle ABO\!=\!\angle OCD\!=\!2\arcsin\frac{a+b}{2}$, 
$\beta\!=\!\angle CBE\!=\!\arccos\left(1-\frac{(b/2)^2 + h^2}{2}\right)$, 
$\gamma\!=\!\angle AEB\!=\!2\arcsin\frac{a}{2}$. 
Calculations show that which of these is larger depends on the ratio $b/a$. For $b<2a$ curving with radii $r$ and $1-r$ is more effective (for $b=a$ the difference is about 13\%), whereas for $b>2a$ curving with radius 1 is more effective  (for $b=5a$ the difference is about 5\%). In our optimizations it always turned out to be best to set $b \approx a$, therefore curving with radii $r$ and $1-r$ is used, giving an increment of area $S_2$.

\paragraph{Pigeonhole principle.}
The order of the $k$-colorable unit-distance graph can be estimated from the known value $\rho_k$ for some fixed tiling based on the pigeonhole principle. The essence of the approach, presented by Pritikin \cite{pri}, is that for a certain number $N$ of vertices of the graph (with an arbitrary arrangement of vertices), one can always choose a tiling position such that all vertices of the graph appear on colored tiles (hence the graph will be $k$-colorable). Technically, this result is ensured by dividing the tiling into many sufficiently small rectangles and by counting the number of shifts that could contain at least one vertex of an arbitrary $N$-vertex graph. Note, that only parallel translations are used here. As a result, the value $\lceil \rho_k-1 \rceil$ (the greatest integer strictly less than $\rho_k$) gives the number of vertices $N$ of an arbitrary unit-distance graph that can be colored in $k$ colors. The value $\lceil \rho_k \rceil$, respectively, gives the lower bound on the order of the $(k+1)$-chromatic graph.

\section{Bounds}

\subsection{$k=1\dots 4$}
For four or fewer colors, the tiling construction of Croft \cite{cro} is used, since tiles of different colors are placed freely in nodes of a hexagonal grid and do not interact with each other (they can be freely translated relative to each other in a certain range).

\begin{figure}[!b]
\centering
\includegraphics[scale=0.34]{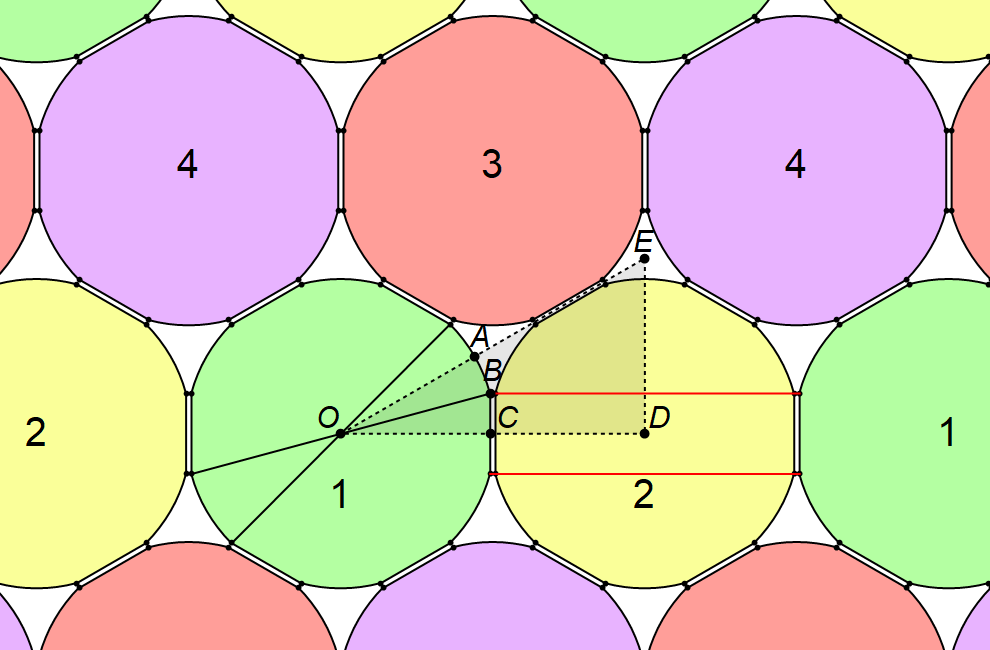}
\caption{Optimal 4-tiling of the plane.}
\label{k4}
\end{figure}

With the optimal shape, the tile has a 12th-order symmetry, so it’s enough to consider the range of angles $\pi/6$, where only two constraints act: one inner and one outer (see Fig.\ref{k4}). 

For $k=1$ the optimal value of the angle $\theta =\angle BOC \approx \pi/12$ is determined by the maximum ratio of the areas of the shapes $OABC$ and $ODE$, which is expressed as $(1-\delta_1)=\left(\pi/6-\theta+\sin(2 \theta)/2\right)/\left((1 + \cos\theta)^2/\sqrt{3}\right)$. Thus $\theta\approx 0.2633155390$, $\delta_1\approx 0.7706352684$, $\rho_1=1.2976307224$.

For $k=1\dots 4$ we get $\delta_k=1-k\,(1-\delta_1)$ and  $\delta_2\approx0.5412705367$, $\delta_3\approx0.3119058051$, $\delta_4\approx0.0825410734$.
Whence $\rho_k=1/\delta_k$ and
$\rho_2\approx1.8475049575$, $\rho_3\approx3.2060961470$, $\rho_4\approx12.1151804529$.

Any departure from this symmetry only worsens the result. It will be extremely surprising if someone can improve this construction.

\subsection{$k=5$}
Surprisingly, no attempt to maximize $\rho_5$ in the case of 5 colors has been published hitherto; the only mention of such work that we are aware of is in Pritikin \cite{pri}, when he wrote that “the results are clumsy and not too impressive” and did not elaborate further. 

\begin{figure}[!b]
\centering
\includegraphics[scale=0.35]{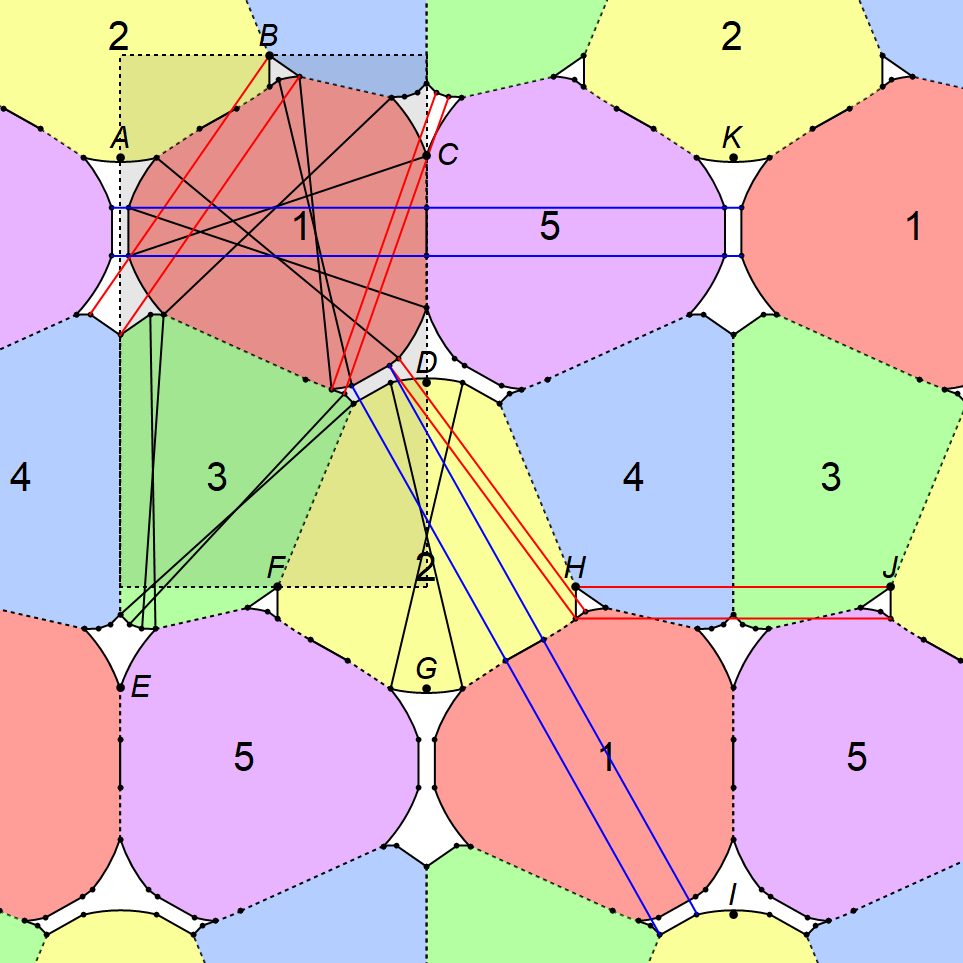}
\caption{Optimal 5-tiling of the plane.}
\label{k5}
\end{figure}

We arrived at the construction shown in Fig.\ref{k5} by starting from Croft’s construction for $k=4$, as follows.
The first attempt is simply to interpose tiles of a 5th color on the Croft tiling so as to cover dumbbell-shaped regions of the void. This gives $\rho_5>16$. 

We can do better if we instead remove the tiles of one color to give the Croft 3-tiling, following which we can place the tiles of the two remaining colors in a symmetric fashion so that they partially overlap each other and also cover larger parts of the original void. This gives $\rho_5>21$ in the optimal arrangement.

To get further, one must change tiles and grid parameters. We did this, progressively complicating the shape of the tiles and introducing outer and inner constraints. Initially, for simplicity, all outer constraints were set by straight strips with unit-width, and all curved sides were described by arcs of radius 1/2 (though one tile may have multiple pairs of arcs, which centers do not necessarily coincide). Interestingly, after removing the restrictions on curving shape, no further improvement emerged, so we omit further details.

Preserving the mirror symmetry of the tiling, we finally came to the construction shown in Fig.\ref{k5} and described by a system of equations with 39 variables. The repeating cell includes 5/2 tiles (two whole ones and one half) and 5/2 voids (one whole and three halves). Taking into account mirror symmetry, one can distinguish 3 types of tiles and 4 types of voids. The number of base points of each type of void is as follows\footnote{After curving, the points $A_1$ and $D_1$ lie in the interior of tiles; they are used only to simplify the calculation of the void area.}: $|\,A\,|=12$, $|\,B\,|=4$, $|\,C\,|=8$, $|\,D\,|=16$. 

To get the area of the voids $S_{\Delta}$ it is convenient to calculate separately the areas of polygons $S$ and curving pairs $S_2$, and take their difference. The area of $n$-gon can be obtained by the coordinates $(x_i,y_i)$ of its corners using triangulation according to the following formula, which allows us to work with negative area increments:
$S=\frac12 \sum_{i=2}^{n-1} \left| \begin{matrix} x_i-x_1 & y_i-y_1 \\ x_{i+1}-x_1 & y_{i+1}-y_1 \end{matrix} \right|$.

The cell area $S_{\Sigma} = x \cdot y$.

Set of constraints\footnote{The points of voids $F$ and $J$, obtained by the reflection of $B$-type void, are numbered counterclockwise.}:
inner (shown black) $C_{1,2}A_{5,4}$, $A_{2,3}D_{9,8}$, $B_{2,3}D_{6,5}$, $D_{2,16}G_{12,2}$, $D_{3,4}E_{6,5}$, $A_{5,6}E_{8,7}$,
outer (red) $B_{1,2}A_{7,8}$, $C_{6,7}D_{4,5}$, $H_{3,4}D_{7,8}$, $H_{1,4}J_{4,1}$,
outer-pair (blue) $D_{6,7}I_{2,3}$, $A_{10,11}K_{3,4}$.

After optimization, an estimate $\rho_5\approx 24.9627819109$ is obtained. All the tiles turn out to be non-convex and to have some wavy sides. Optimally rigid tile boundaries are shown solid; non-rigid ones are dashed. The numbers of corners of tiles with wavy sides are: 14 for tiles of colors 1 and 5, 10 for tiles 2, 12 for tiles 3 and 4. The minimal number of tile corners after dividing wavy edges into line segments or arcs depends on the form of non-rigid boundaries and can be 15, 14, and 12, respectively.

The fraction of the plane occupied by voids in our construction is $\delta_5 \approx 0.040059637727$. It is a pity that there is not enough trifle to a beautiful mark $\delta_5=4\%$. Maybe someone see where else could we push?

\subsection{$k=6$}
In contrast to $k=5$, prior work on the case $k=6$ has been quite extensive, and we use it as our starting-point. The first tiling of nearly all of the plane with six colors was obtained by Pegg \cite{soi}, with $\rho_6\approx 303$. In his construction, all tiles are heptagons and the voids are squares. Pritikin \cite{pri} explored a generalization of Pegg’s construction where the voids are rhombi and obtained a large improvement, to $\rho_6\approx 6197$. Here we further refine Pritikin’s construction, obtaining a larger value for $\rho_6$.

\begin{figure}[H]
\centering
{
\centering
\begin{tabular}{cccc}
    \multicolumn{4}{c}{\includegraphics[scale=0.305]{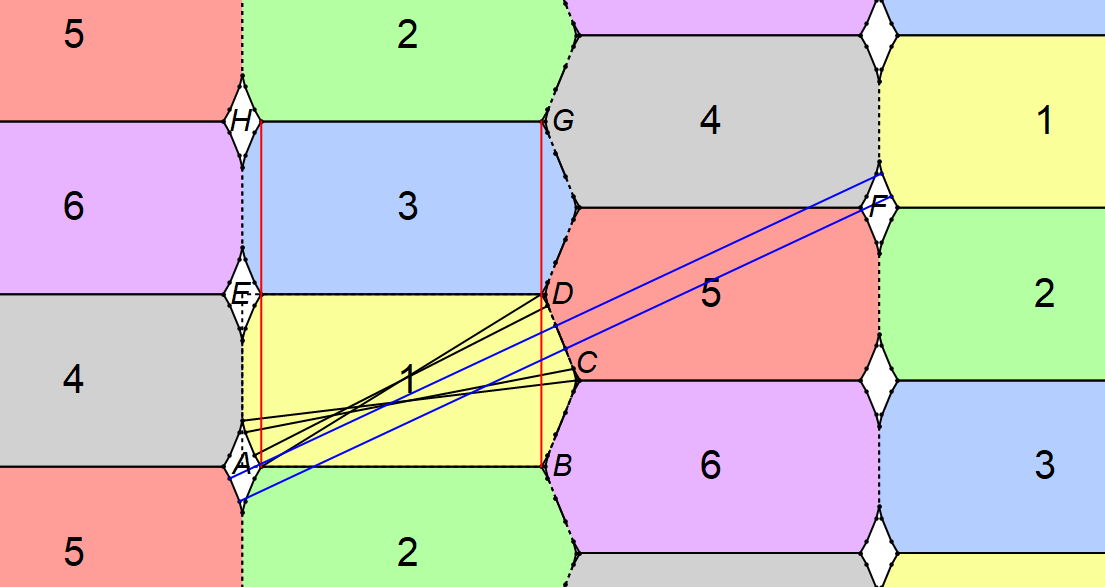}} \\ [8pt]
    \includegraphics[scale=0.3]{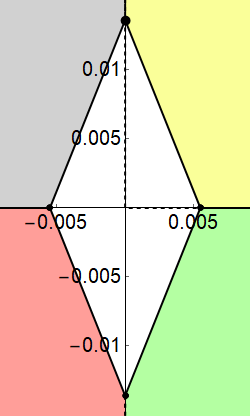} & \includegraphics[scale=0.3]{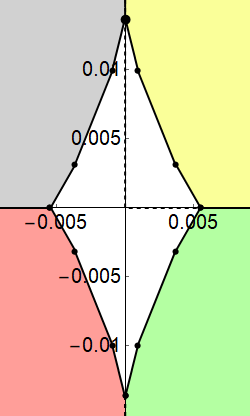} & \includegraphics[scale=0.3]{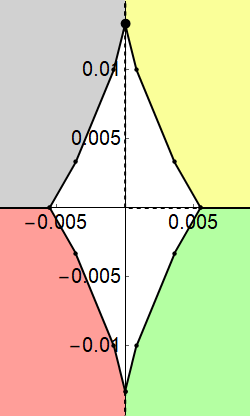} & \includegraphics[scale=0.3]{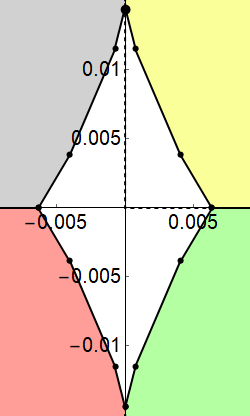} \\ [8pt]
    \includegraphics[scale=0.3]{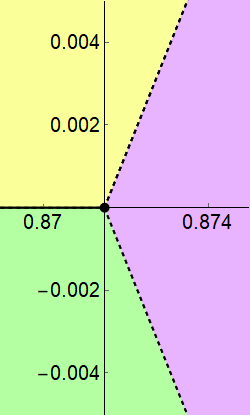} & \includegraphics[scale=0.3]{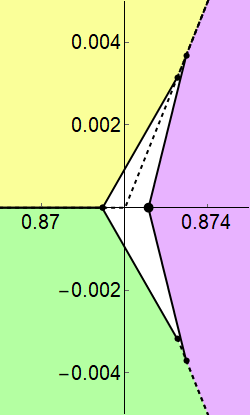} & \includegraphics[scale=0.3]{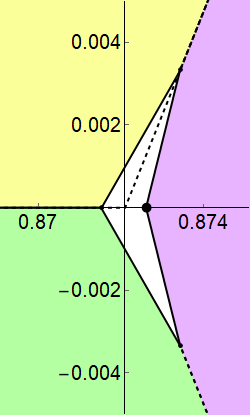} & \includegraphics[scale=0.3]{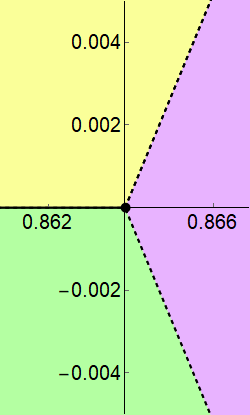} 
\end{tabular} \par
}
\caption{Optimal 6-tiling of the plane. Above is a general view of the optimal tiling, with the dimensions of all voids increased by 10 times. Below are the real sizes of voids of type $A$ (middle row) and $B$ (lower row) for different tile structures, from left to right: 1) Pritikin’s 7-gons, $\rho_6\approx6197$, 2) convex 15-gons, $\rho_6\approx6906$, 3) non-convex curved 17-gons, $\rho_6\approx6992$, 4) convex curved 11-gons, $\rho_6\approx6043$ (7-colorable tiling).}
\label{k6cons}
\end{figure}

We start by tiling the plane with pentagons (see Fig.\ref{k6cons}). In the construction of Pritikin, the heptagonal tile fits entirely inside the pentagon $ABCDE$ with coordinates \{$(0, 0)$, $(x-z, 0)$, $(x+z, y/2)$, $(x-z, y)$, $(0, y)$\}. The angle $\angle BAE$ is right, the slope of the edges $BC$ and $CD$ is defined by the macro parameters $\{y, z\}$, the repeating cell has an area $S_{\Sigma} = x \cdot y$.

In this tiling, two types of vertices can be distinguished: $A$ and $B$, with degree 4 and 3 respectively. To obtain a partial 6-tiling, rhombus-shaped voids are introduced at vertices of type $A$. (The dimensions of these voids are so small that they are almost indistinguishable in the figure, if we follow the scale. To see all the most interesting, we need a microscope.) We denote the coordinates of corners of rhombus $A$ as \{$(0, l)$, $(r, 0)$, $(0, -l)$, $(-r, 0)$\}.

The original construction uses heptagons with straight sides, giving the following constraints: inner $A_1C$, $A_2D$, outer $A_2BGH_2$, outer-ext $A_{3,4}DC$; $S_{\Delta} = l\cdot r$. Optimum result\footnote{Pritikin gave a slightly suboptimal estimate of 6197.08. He specified exact decimal values for the dimensions of the rhombi; since his main focus was the lower bound on the order of a 6-colorable unit-distance graph, we may assume that he chose to specify only enough significant digits to achieve the maximum integer part for $\rho_6$.} $\rho_6 \approx 6197.5793083297$.

To improve this, we first note that the sides $A_3A_4$ and $CD$ are not parallel; therefore, it is possible to improve the result by allowing a wavy side $CD$, and by making its middle part parallel to the side $A_3A_4$. (As a result, the tiles cease to be convex, and the outer-ext constraint becomes an outer-pair one.)

Next, one can add arrowhead-shaped voids at the corners of the inclined sides instead of the vertices of type $B$. 

It turns out that this allows a reduction in the size of the rhombi that exceeds the area of the arrowheads themselves, giving an 
increase in $\rho_6$.

Next, one can introduce additional corners on the rhombi, converting them into dodecagonal voids. This removes the inefficiency associated with requiring each corner of a rhombus to obey two constraints (an inner one and an outer-pair one).

The resulting tile is a 15-gon
$A_{1,2,3,4}B_{4,5}C_{6,1,2}D_{3,4}E_{4,5,6,7}$ = \{$(0, l)$, $(p, m)$, $(q, n)$, $(r, 0)$, $(x-z-s, 0)$, $(x-z+4vz/y, v)$, $(x+z-4wz/y, y/2-w)$, $(x+z-t, y/2)$, $(x+z-4wz/y, y/2+w)$, $(x-z+4vz/y, y-v)$, $(x-z-s, y)$, $(r, y)$, $(q, y-n)$, $(p, y-m)$, $(0, y-l)$\}, with  $S_{\Delta}=lp+nr+mq-np+vs+wt$, $S_{\Sigma}=xy$, and the following set of constraints:
inner $A_{1,2}C_{2,1}$, $A_{3,4}D_{4,3}$,
outer $A_4B_4G_4H_4$,
outer-pair $A_{8,9}F_{2,3}$.
In total, 13 variables are used: $\{l, m, n, p, q, r, s, t, v, w, x, y, z\}$.

Finally, one can curve four pairs of sides of each tile, which slightly improves the result. These are $A_{1,2}C_{2,1}$, $A_{3,4}D_{4,3}$, $E_{4,5}B_{5,4}$ and $E_{6,7}C_{1,6}$. 

With all these refinements implemented, we obtain an optimum value
of $\rho_6 \approx 6992.1655504123$. The coordinates of the tile’s base points are given by the following parameter values:
$x \approx 0.9220971880$, $y=0.5$, $z \approx 0.0500073975$,
$l \approx 0.0133435475$, $m \approx 0.0099987220$,
$n=v=w \approx 0.0033448254$,
$p \approx 0.0008269367$, $q \approx 0.0035745405$, $r \approx 0.0054885277$, 
$s \approx 0.0005758591$, $t \approx 0.0005111915$.
(It should be mentioned that we did not initially require these quadrilaterals to be rectangles or even trapezoids, but that rectangles turned out to be the optimal shape. Similarly, no improvement was obtained by allowing $A_{8,9}F_{2,3}$ to be a parallelogram.)

Table \ref{tpri} shows the optimal results when various subsets of the optimizations are included. The greatest effect is achieved by adding type $B$ voids (about +6.4\%), followed by the increase in the number of corners of type $A$ voids (+5.0\%), wavy edges (+1.6\%) and curving (+0.1\%). The overall effect is an increase in $\rho_6$ of about 12.8\% compared to the construction of Pritikin.

\begin{table}[!t]
{
\caption{Optimal values of $\rho_6$ for various tiling options.}
\label{tpri}
\smallskip

{
\centering
\small
\begin{tabular}{@{}c<{\!\!\!}|>{}c<{}*{2}{|@{}>{\,}c<{\,}@{}}*{4}{|>{\!\!}c<{\!\!}}}  

\hline
 proper & \multicolumn{3}{c|}{sides} & \multicolumn{2}{c|}{non-wavy} & \multicolumn{2}{c}{wavy} \\
\cline{2-8}
 colors & tile & rhmb & arrow & straight & curved & straight & curved \\
\hline
\hline
7 &  7 &  4 & 0 & 5681.489884 &             & 5780.207842 &             \\
  & 11 & 12 & 0 & 5942.027491 & 5943.447950 & 6041.630152 & 6043.112468 \\
\hline
8 &  7 &  4 & 0 & 6197.579308 &             & 6294.390841 &             \\
  & 11 & 12 & 0 & 6510.512176 & 6512.380146 & 6604.352015 & 6606.168990 \\
\hline
9 & 11 &  4 & 6 & 6596.061280 & 6597.901152 & 6684.123145 & 6685.917214 \\
  & 15 & 12 & 6 & 6899.423068 & 6906.361529 & 6985.378399 & 6992.165550 \\
\hline
\end{tabular}

}
}
\end{table}

Note that in the fully optimized version $v=w$
, so that the arrows become quadrilaterals instead of hexagons. As in the case $k=5$, the optimal curving radius is always $1/2$ (i.e., all relevant trapezoids degenerate into rectangles). The minimal number of tile corners if each wavy edge is replaced by three segments is 17 (or 19 in the case of $v\ne w$). The fraction of the plane occupied by voids is $\delta_6 \approx 0.000143017209$.

For completeness, we finally note that one can dispense with voids of type $A$ entirely. If we use only voids of type $B$, we get $\rho_6 \approx 282.1748985451$. In this case, the void boundaries are formed by arcs of radius 1.

\subsection{$k=7$}

Proper 7-tilings of the plane have been known since the 1950s, and an infinite number of such tilings exist, so values for $k$ greater than 6 would seem uninteresting for present purposes. However, it is worth exploring whether our constructions might facilitate the search for a proper 6-tiling, whose existence remains an open problem \cite{soi}. To do so, we can define quantitative metrics of proper 7-tilings that tend to zero in the limit of a proper 6-tiling, so as to estimate how far we are from tiling the entire plane using 6 colors only.

One such metric is the minimum value, among the seven colors, of the fraction of the plane that tiles of that color occupy. In the original construction of Pegg \cite{soi}, the voids are purposely spaced such that they can all be assigned the same color. In Prikitin’s refinement \cite{pri} this is not so, and alternate voids (of type $A$) must be assigned different colors, so 8 colors are needed in total; and our optimal tiling (with curved 17-gons and two types of voids $A$ and $B$) is 9-colorable (guess how). We can obtain a proper 7-tiling by eliminating the arrowhead-shaped voids and slightly increasing the width of the tiles ($A_1H_7$ is now required to be at least 1), but retaining our other optimizations. In this case the optimal parameter values are:
$x \approx 0.9169683902$, $y \approx 0.5143908364$, $z \approx 0.0531068312$, $l \approx 0.0143908364$, $m \approx 0.0115262526$, $n \approx 0.0038562709$, $p \approx 0.0007214876$, $q \approx 0.0040043282$, $r\approx 0.0063056906$, $s=t=v=w=0$,
leading to\footnote{Here we are actually consider partial $6$-tiling, so the notation $\rho_7$ would be incorrect.} $\rho_6 \approx 6043.11246787913$.
This improves Pegg's result $\rho_6 \approx 303.0962048019$ \cite{soi}\footnote{To pre-empt any confusion we note that the value given in \cite{soi} for Pegg’s construction is the ratio of tile area to void area, so $\rho_6-1$.}. The fraction of the seventh color in our case is $\delta_6 \approx 0.000165477642$.

Another possible metric is the minimum ratio, among the seven colors, of tiles of at color to the total number of tiles. 

\begin{figure}[H]
{
\centering
\begin{tabular}{cc}
    $a$) & \includegraphics[scale=0.47]{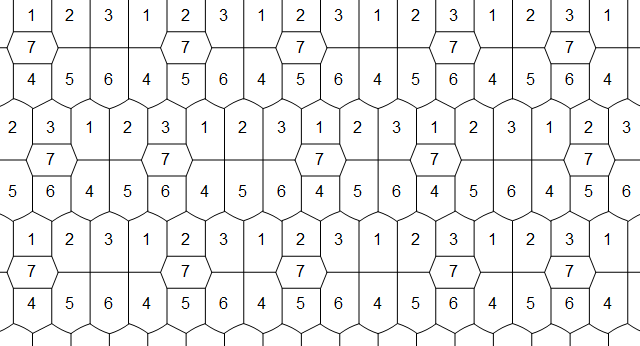} \\ [10pt]
    $b$) & \includegraphics[scale=0.47]{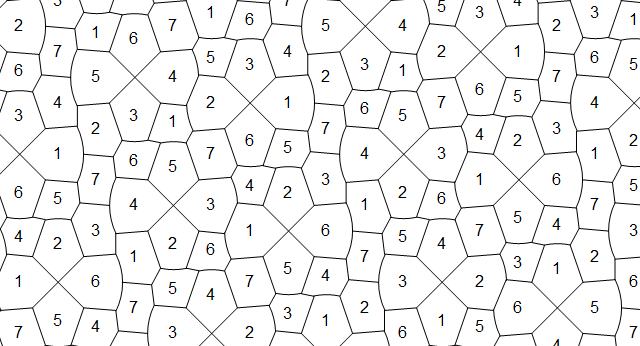} \\ [10pt]
    $c$) & \includegraphics[scale=0.47]{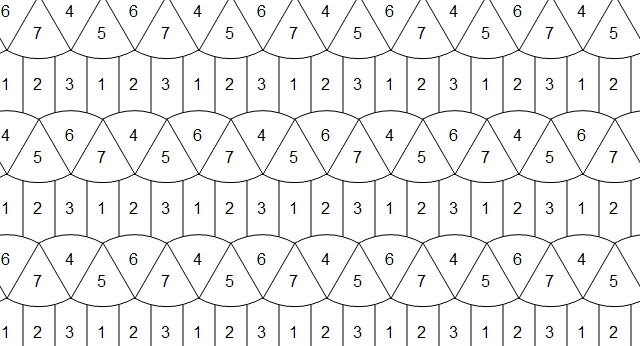} 
\end{tabular} \par
}
    \caption{Examples of 7-tilings with reduced occurence of 7th color.}
    \label{k7}
\end{figure}

Quite unexpectedly, the current estimates here are much worse. We are unaware of any way to obtain a Siamese-free tiling with the relative proportion of tiles of some color being less than 1/8.

Fig.\ref{k7} shows some 7-tilings, such that tiles of the 7th color have occurrence frequency 1/8: 
a) based on the tiling of Pritikin-Pegg, the ratio of the number of tiles is 7:7:7:7:7:7:6,
b) obtained on the basis of the hexagonal tiling of the plane, also 7:7:7:7:7:7:6,
c) with Reuleaux triangles used to separate rows of pentagons, 4:4:4:3:3:3:3. 

This result (or lack of it) allows us to take a different look at the gap between 6- and 7-colorings. If we start from the partial 6-tiling described in the previous section, we may get the feeling that we are very close to a proper 6-tiling of the plane. However, if we pay attention not to the occupied area, but to the relative number of tiles of a given color (which is 1/7 for the classical honeycomb-based 7-tiling), it seems that we are very far from at goal.

\subsection{Observations}
For calculations, we used the package {\tt Mathemetica 10}, in particular, the function {\tt FindMaximum}. 
It was found that constraints should be set in the form of equalities. Even if inequalities are partially used, the convergence of the solution is significantly impaired. It is also important to choose a good initial approximation (which can be obtained in the previous step using a simpler tiling construction).

The curving of tile sides was carried out jointly with the optimization of the position of tile corners. If the shape of the polygons is optimized first, and then the sides are curved, the final result is worse.

In our investigations, optimal tilings for 5 and 6 colors are characterized by a number of features: 
a) vertices in which three or more colored tiles converge are disadvantageous (a void is superior to any such vertex);
b) tiles become non-convex (i.e., with wavy edges);
c) pairs of parallel sides of the same length are formed (so that trapezoids with opposite edges that can be curved always degenerate into rectangles);
d) the curving of pairs of opposite tile sides is optimized with a radius of 1/2.

\section{Future work}
Here we have restricted ourselves to the case of highly symmetrical tilings with only a few tiles in the repeating cell, and in the case of $k=6$ to tilings in which all tiles are congruent. There is a chance that using more asymmetric tilings, it will be possible to increase the fraction of the colored plane. But a thorough search for such tilings requires the development of more advanced analysis methods. When trying to improve the results presented here, it may be most efficient to incorporate observed features listed above as assumptions at the outset.

It is very doubtful that Siamese tiles can help in finding a proper 6-tiling, since in this case it is necessary to choose the appropriate partial 5-tiling, such that all voids will satisfy the unit-distance requirement. However, it should be recognized that in some special circumstances Siamese pairs of tiles may be fruitful. For example, our proper 7-tiling derived from partial 6-tiling contains them.

Table \ref{tspec} shows the known bounds of the order of the minimal $k$-chromatic unit-distance graph. The lower bounds are obtained on the basis of tiling by appealing to Pritikin’s argument; the upper bounds arise from construction of the corresponding graph. One may wonder whether it is possible to predict the size of a 6-chromatic graph (if it exists)?

\begin{table}[!t]
{
\caption{Bounds on the order of unit-distance graphs in the plane.}
\label{tspec}
\smallskip

{
\centering
\small
\begin{tabular}{l|c|c|c|c|c|c|c}
\hline
number of colors, k              & 1 & 2 & 3 &  4 &  5 &    6 & 7 \\
\hline
k-colorable graph by pigeonhole  & 1 & 1 & 3 & 12 & 24 & 6992 & $\infty$ \\
k-chromatic graph by pigeonhole  &   & 2 & 2 &  4 & 13 &   25 & 6993 \\
smallest known k-chromatic graph & 1 & 2 & 3 &  7 & 509 &   &   \\	\hline
\end{tabular}

}
}
\end{table}

It can be noted that even for small $k$, the bounds obtained using the pigeonhole principle are very weak and quickly diverge from the exact (for $k\le 4$) values with increasing $k$. We can only expect this trend to worsen with larger $k$, even if we assume that there is a much more compact graph for $k=5$ \cite{par}, which is hard to believe.

Here are some obvious reasons showing why pigeonhole-based lower bounds are a bad guide. Important details that it does not take into account are: a) rotations when covering the graph by tiling (in addition to parallel translations), b) the shape and location of tiling voids, c) the uneven Euclidean density of the vertices of finite graphs (higher density nearer the center), d) the degree of graph vertices, which should be at least $k$ (for an arbitrary set of vertices it tends to zero).

Taking these factors into account should increase the estimate of the number of vertices. Perhaps better approximations 
can be obtained by considering some combination of tiling the disk and the entire plane. For example, one can introduce a weight function that penalises the appearance of voids by an amount that decreases with distance from the center. In addition, if we recall the history of the discovery of 5-chromatic graphs, then after obtaining the first such graph \cite{grey}, the minimal order was reduced by about 40 times, or even more. Therefore, when searching for 6-chromatic graphs, a significant margin must be laid.

As long as we do not have more adequate bounds, we can focus on existing ones with some correction. Here everyone is free to offer their own version of data extrapolation. We want to think that “the soil on the moon is solid”\footnote{
When the first steps of space exploration were taken and no one knew what to expect when landing on the moon, General Designer S.P. Korolev took responsibility upon himself and issued the following statement: “The landing of the lunar spaceship should be expected on a fairly hard surface, such as pumice... 28/X/64, Korolev”
} and there exist 6-chromatic unit-distance graphs with a number of vertices less than 1 million, which would give hope to reach them with SAT methods. So it's time to take off.

\section{Acknowledgements}
Alexander Soifer did the beginning and took responsibility for the translation of the most important part of the article. Aubrey de Grey did the rest (I just got in the way). Why am I in the title instead of them, and why alone, it's a mystery for me.

\end{document}